\documentclass[12pt]{amsart}
\usepackage{amssymb,amsmath,amscd,graphicx,
latexsym,amsthm}
\usepackage{amssymb,latexsym,eufrak,amsmath,amscd,graphicx,array}
\usepackage[mathscr]{eucal}
  \usepackage[all]{xy}

\usepackage[all]{xy} 

\renewcommand{\P}{ \ensuremath{\mathbb{P}}}

\usepackage{layout}

\theoremstyle{plain}
\newtheorem*{theoremA}{Theorem A}
\newtheorem{theorem}{Theorem}[section]
\newtheorem{lemma}[theorem]{Lemma}
\newtheorem{Proposition}[theorem]{Proposition}

\theoremstyle{definition}
\newtheorem{definition}[theorem]{Definition}
\newtheorem*{definitionA}{Definition}
\newtheorem{remark}[theorem]{Remark}
\newtheorem{example}[theorem]{Example}

\newcommand{\noi}{\noindent}
\begin{document}
\title{Regularity of smooth curves in biprojective spaces}

\author{Victor Lozovanu} 
\address{Department of Mathematics \\ University of Michigan \\ Ann Arbor 
\\ MI 48109 \\ USA}
\email{vicloz@umich.edu}

\date{16 October 2008}

\setcounter{section}{-1}
\begin{abstract}
Maclagan and Smith \cite{MaclaganSmith} developed a multigraded version of Castelnuovo-Mumford regularity. Based on their definition we will prove in this paper that for a smooth curve $C\subseteq \P^a\times\P^b$ $(a, b\geq 2)$ of bidegree $(d_1,d_2)$ with nondegenerate birational projections the ideal sheaf $\mathcal{I}_{C|\P^a\times\P^b}$ is $(d_2-b+1,d_1-a+1)$-regular. We also give an example showing that in some cases this bound is the best possible.
\end{abstract}

\maketitle

\section*{Introduction}
Let $\mathcal{F}$ be a coherent sheaf on $\P^r$. Recall (Lecture 14, \cite{Mumford66}) that $\mathcal{F}$ is $m$-regular if 
\[
H^i(\P^r,\mathcal{F}\otimes \mathcal{O}_{\P^r}(m-i))=0, \forall i>0
\]
One then defines: reg$(\mathcal{F})=$min$\{m\in \mathbb{Z} :\mathcal{F} \text{ is $m$-regular} \}$. This invariant bounds the algebraic complexity of a coherent sheaf, and for this reason has been the focus of considerable activity, e.g. \cite{C}, \cite{Mumford70}, \cite{GP1}, \cite{GP2}, \cite{GLP}, \cite{L}, \cite{EL},\cite{K98}, \cite{K00}.

In recent years , it became of great interest to extend the definition of regularity to the multigraded case and study it's behaviour, e.g. \cite{MaclaganSmith}, \cite{HT}, \cite{HW}, \cite{ST}, \cite{STW}, \cite{H}, \cite{HS}. Maclagan and Smith \cite{MaclaganSmith}, working on toric varieties, introduced a notion of multigraded regularity, that on $\P^a\times\P^b$ reduces to:
\begin{definitionA}
We say that a coherent sheaf $\mathcal{F}$ on $\P^a\times\P^b$ is $(m,n)$-regular if:
\[
H^i(\P^a\times\P^b,\mathcal{F}\otimes \mathcal{O}_{\P^a\times\P^b}(m-u,n-v))=0
\]
for all $i>0$, and all  $(u,v)\in\mathbb{N}^2$ with $u+v=i$.
\end{definitionA}
\noi Let reg$(\mathcal{F})=\{ (m,n)\in \mathbb{Z}^2 : \mathcal{F} \text{ is $(m,n)$-regular }\}$. Although in contrast with the usual regularity, reg$(\mathcal{F})$ is a set rather than a single number, they still share common properties. For example, if $\mathcal{F}$ is $(0,0)$-regular sheaf then it is globally generated \cite{MaclaganSmith}, \cite{HSS}. 

In the classical case an important problem is to find bounds on Castelnuovo-Mumford regularity for subvarieties of projective space. As an example  Gruson, Lazarsfeld and Peskine \cite{GLP} proved that when $C\subseteq\P^r$ $(r\geq 2)$ is an irreducible, nondegenerate curve of degree $d$ then $\text{reg}(C):=\text{reg}(\mathcal{I}_{C|\P^r})\leq d-r+2$. Now it is natural to ask if we can give bounds on the multigraded regularity of a smooth curve $C\subseteq\P^a\times\P^b$ depending on its bidegree. The purpose of this paper is to prove the following theorem:
\begin{theoremA} \label{Theorem A}
Let $C\subseteq\P^a\times\P^b$ $(a, b\geq 2)$ be a smooth curve of bidegree $(d_1,d_2)$ with nondegenerate birational projections then the ideal sheaf $\mathcal{I}_{C|\P^a\times\P^b}$ is $(d_2-b+1,d_1-a+1)$-regular.
\end{theoremA}
\noi As a corollary, Theorem A together with (Theorem $1.4$ \cite{MaclaganSmith}) imply the inclusion:
\[
(d_2-b+1,d_1-a+1)+\mathbb{N}^2 \ \subseteq \text{reg}(C) \ := \text{reg}(\mathcal{I}_{C|\P^a\times\P^b})
\] 
The proof of the theorem uses generic projections and vector bundles techniques developed by Gruson, Peskine(\cite{GP1}, \cite{GP2}) and Lazarsfeld (\cite{L}) in the classical case. The minor change is that instead of projecting to $\P^2$ we choose projections to $\P^1\times\P^1$. We start the paper by proving that whenever $a\neq b$ or $r:=a=b$ and the curve $C$ is not included in the graph of an automorphism of $\P^r$, then there are plenty of projections to $\P^1\times\P^1$, with the image of $C$ having "nice" singularities. In section two these projections will play an important role in establishing the bound from Theorem A. The remaining case when the curve is included in the graph of an automorphism of $\P^r$ will be discussed in section three. There we will use some results from \cite{GLP} to show that the same bound works.

The paper ends in \S 4 with an example of a rational curve $C\subseteq\P^2\times\P^2$ of bidegree $(3,3)$. This curve has the property that $(2,2) + \mathbb{N}^2 = \text{reg}(C)$, showing in this case that the bounds we have in Theorem A are the best possible.

\subsection*{Acknowledgments.}  I would like to thank my advisor, Rob Lazarsfeld, for suggesting me this problem, for his support and many useful comments. I am also grateful to Igor Dolgachev and Mihai Fuger for having the patience to listen to the final draft of the paper; to Jessica Sidman for making the point that I should find an example which supports the maximality of my bounds and to Greg G. Smith for encouragements.  I would also like to thank Fidel Jimenez and Parsa Bakhtary for useful discussions and encouragements.

\section{Definitions and notations}

\noi (0.1). We work throughout over the complex numbers.

\noi (0.2). Unless otherwise stated a \emph{curve} is a \emph{smooth}, \emph{irreducible} projective variety of dimension one. Let $C\subseteq \P^a\times\P^b$ $(a, b\geq 2)$ be a curve and denote by $p_1$ and $p_2$ the projections to each factor. We say $C$ has \emph{nondegenerate birational projections} if $p_1|_C$ and $p_2|_C$ are birational morphisms and have nondegenerate images, call them $C_1$ and $C_2$ respectively. Write $L_1:=\mathcal{O}_{\P^a\times\P^b}(1,0)\otimes \mathcal{O}_C$ and $L_2:=\mathcal{O}_{\P^a\times\P^b}(0,1)\otimes \mathcal{O}_C$. These are two line bundles on $C$ and let $W_i\subseteq H^0(C,L_i)$ be the linear subsystem defining the restriction map $p_i|_C$, where $i=1,2$. The curve $C$ is said to be of bidegree $(d_1,d_2)$, if $d_1=\text{deg}_C(L_1)$ and $d_2=\text{deg}_C(L_2)$.

\noi (0.3). Usually $\Lambda_a$ and $\Lambda_b$ will mean codimension two planes in $\P^a$ and $\P^b$ respectively. They  define the projection maps: $\pi_a:\P^a\dashrightarrow\P^1$ and $\pi_b:\P^b\dashrightarrow\P^1$. If $\Lambda_a \cap C_1=\{\varnothing\}$ and $\Lambda_b\cap C_2=\{\varnothing\}$ then we have the diagram: 
\[
\xymatrix{
C \ar[d]^{f_{a,b}} & \subseteq & \P^a\setminus \{ \Lambda_a\}\times\P^b\setminus \{ \Lambda_b\} \ar[d]^{\pi_a\times\pi_b}\\
\bar{C} & \subseteq & \P^1\times\P^1\\}
\]
where $f_{a,b}$ is the restriction to $C$ of $\pi_a\times\pi_b$. In this case we introduce the definition:
\begin{definition}\label{Definition 2}
We will say that the pair $(\Lambda_a,\Lambda_b)$ defines a \textit{good projection} for $C$ if $f_{a,b}$ is a birational morphism with fibers of length at most two and the differential map of $f_{a,b}$ is injective for all $z\in C$.
\end{definition}
\noi (0.4) Usually $H_a$ and $H_b$ will mean hyperplanes in $\P^a$ and $\P^b$ respectively. They define the following finite sets:
\[
(C.H_b):=p_1(\P^a\times H_b\cap C)\subseteq\P^a\text{ and }(H_a.C):=p_2(H_a\times\P^b\cap C)\subseteq\P^b
\]
\noi Consider two subsets $F_1,F_2\subseteq\P^a$ then we will denote by $\overline{F_1F_2}$ the union of all lines connecting one point from $F_1$ and another one from $F_2$.

\section{Existence of good projections}
In this section we prove the existence of good projections. Specifically our goal is to establish the following theorem:
\begin{theorem} \label{Theorem 1}
Let $C\subseteq \P^a\times\P^b$ be a smooth curve with nondegenerate birational projections. Suppose that either $a\neq b$ or $r:=a=b$ and the curve $C$ is not included in the graph of an automorphism of $\P^r$. Then $C$ has good projections to $\P^1\times\P^1$.
\end{theorem}
\begin{remark} \label{Remark 1}
In search for good projections it is necessary that a general center $\Lambda_a\in\text{Grass}(a-2,\P^a)$ is not contained in a hyperplane $H_a\subseteq\P^a$, where the map $p_2$ projects at least two points of the set $H_a\times\P^b\cap C$ to the same one. A hyperplane like this has two points $(x_1,y_1),(x_2,y_2)\in H_a\times\P^b\cap C$ with $y_1=y_2$ and $x_1\neq x_2$ and as a result $H_a$ should contain the line $\overline{x_1x_2}$.  We assumed that $p_2|_C$ is birational to its image so there exist only finitely many pairs of points on $C$ having the same image under $p_2$. This implies that the family of these hyperplanes is of codimension at least two. As a consequence the dimension of those $\Lambda_a\in\text{Grass}(a-2,\P^a)$ contained in these hyperplanes is at most $a-2+a-1=2a-3$. But dim$(\text{Grass}(a-2,\P^a))=2a-2$ and the assertion follows immediately.
\end{remark}
The proof of Theorem \ref{Theorem 1} uses the idea of "uniform position principle" developed by Harris, e.g. Chapter III of \cite{ACGH}. Specifically  we have the following lemma:
\begin{lemma} \label{Lemma1}
If $\Lambda_a\subseteq \P^a$ is a general codimension two plane then one of the following two situations must happen:
\begin{enumerate}
\item[$(1_a)$] For all hyperplanes $H_a$ containing $\Lambda_a$, the set $(H_a.C)$ does not span $\P^b$.
\item[$(2_a)$] For a general hyperplane $H_a$ containing $\Lambda_a$, any $b+1$ points from $(H_a.C)$ span $\P^b$.
\end{enumerate}
\end{lemma}
\begin{proof} 
The curve $C$ is the desingularization of $C_1$, then (Theorem 1.1 \cite{PS}) implies that the projection map $\pi_a|_C$ defined by a general codimension two plane $\Lambda_a\subseteq\P^a$ has the monodromy the full symmetric group. If we set
\[
U=\P^1\setminus\{\textrm{Branch points of }\pi_a|_C\} \text{ and }V=\pi^{-1}_a|_C(U)
\]
this says that $\forall y\in U$ every two points in the fiber $\pi^{-1}_a|_C(y)$ can be connected by a path in $V$ lifted from a loop in $U$ based at $y$. Now construct the following incidence correspondence:
\[
I_a(b+1)\subseteq V\times....\times V\times U
\]
consisting of those tuples $(q_1,...,q_{b+1},y)$, where the points $q_1,...,q_{b+1}$ are distinct and contained in the fiber $\pi^{-1}_a|_C(y)$. As the monodromy is the full symmetric group, $I_a(b+1)$ is connected. At the same time the projection map $I_a(b+1)\rightarrow U$ is a covering space and $U$ is irreducible, therefore $I_a(b+1)$ is an irreducible variety of dimension one.

Let's define the following closed subvariety:
\begin{eqnarray}
J_a(b+1) & = &\{(q_1,...,q_{b+1},y)\in I_a(b+1): \nonumber\\
&& p_2(q_1),...,p_2(q_{b+1}) \text{ don't span } \P^b\}\nonumber
\end{eqnarray}
But $I_a(b+1)$ is irreducible, hence either $J_a(b+1)=I_a(b+1)$ or dim$(J_a(b+1))=0$. Bearing in mind Remark \ref{Remark 1} we have that the first case corresponds to $(1_a)$ and the second one is equivalent with $(2_a)$.
\end{proof}
\begin{lemma} \label{Lemma2}
Let $\Lambda_a\subseteq \P^a\setminus C_1$ be a general codimension two plane in the first factor. If condition $(2_a)$ from Lemma \ref{Lemma1} is satisfied then for a general codimension two plane $\Lambda_b\subseteq \P^b$ the pair $(\Lambda_a,\Lambda_b)$ defines a good projection for $C$.
\end{lemma}
\begin{proof}
Since $\Lambda_a$ is general, Remark \ref{Remark 1} tells us that $\forall x\in\P^1$ both $\overline{x\Lambda_a}\times\P^b\cap C$ and $(\overline{x\Lambda_a}.C)$ have the same cardinality. With this in hand we can start the proof and first we need for a general choice of $\Lambda_b$, the map $f_{a,b}$ to have the fibers of length at most two. The map $f_{a,b}$ has a fiber of length at least three if there exist $x\in\P^1$ and a hyperplane $H_b$ passing through $\Lambda_b$, which contains at least three points from the set $(\overline{x\Lambda_a}.C)$. Condition $(2_a)$ guarantees the existence of an open set $U\subseteq\P^1$, where $\forall x\in U$ any three points in $(\overline{x\Lambda_a}.C)$ span a plane in $\P^b$. Hence the family of those hyperplanes $H_b$ which for some $x\in U$ contain at least three points from the set $(\overline{x\Lambda_a}.C)$ is of codimension two. Simultaneously $\Lambda_b$ should not be included in a hyperplane which for some $x\in\P^1\setminus U$ contains at least two points from $(\overline{x\Lambda_a}.C)$. Now beairing in mind the ideas from Remark \ref{Remark 1}, the assertion follows immediately.

Next, we need for a general $\Lambda_b$ the map $f_{a,b}$ to be birational to its image. The map is not birational if given a general hyperplane $H_b$ containing $\Lambda_b$, there exists $x\in \P^1$ such that $H_b$ contains at least two points from the set $(\overline{x\Lambda_a}.C)$, hence forcing $\Lambda_b$ to intersect the line connecting these two points. But the union of all lines connecting two points of $(\overline{x\Lambda_a}.C)$ when $x\in\P^1$ is of dimension two; so we have the birationality condition as well.

Lastly we need the differential map of $f_{a,b}$ to be injective for all $z\in C$. Then there should not exist a hyperplane $H_b$ passing through $\Lambda_b$ and a point $x\in\P^1$ such that $\overline{x\Lambda_a}\times H_b$ contains the tangent direction at some point on $C$. The fact that $\Lambda_a$ is general implies that for only finitely many $x\in\P^1$ we have that $\overline{x\Lambda_a}\times\P^b$  is tangent to $C$ at some point. Therefore a general $\Lambda_b$ does not intersect the projection to $\P^b$ of the tangent direction at these points so the differential map of $f_{a,b}$ is injective for all $z\in C$.
\end{proof}
Lemma \ref{Lemma1} and Lemma \ref{Lemma2} imply that we might not be able to obtain good projections, only if condition $(1_a)$ in Lemma \ref{Lemma1} is satisfied for a general codimension two plane $\Lambda_a\subseteq\P^a$. This means that for a general hyperplane $H_a\subseteq\P^a$, the set $(H_a.C)$ does not span $\P^b$, and as this condition is closed then it is true for all hyperplanes. Now if we start the choice with a codimension two plane in $\P^b$, we deduce that the only case when we might not be able to produce good projections is when the curve $C$ satisfies the property:
\begin{center}
\begin{tabular}{cm{5.5in}}
$(*)$ & For all hyperplanes $H_a\subseteq\P^a$ and $H_b\subseteq\P^b$, the finite sets $(H_a.C)$ and $(C.H_b)$ don't span $\P^b$ and $\P^a$ respectively.
\end{tabular}
\end{center}
\noi The last fact needed to prove Theorem \ref{Theorem 1} is the following lemma:
\begin{lemma}\label{Lemma3}
Let $C\subseteq\P^a\times\P^b$ be a curve which satisfies $(*)$, then for a general hyperplane $H_a\subseteq\P^a$ we have:
\begin{enumerate}
\item[$(i)$] The set $H_a\cap C_1$ spans the hyperplane $H_a$.
\item[$(ii)$] The set of points $(H_a.C)$ spans a unique hyperplane $H_b$ in $\P^b$.
\item[$(iii)$] $H_a\times\P^b\cap C=\P^a\times H_b\cap C$.
\end{enumerate}
\end{lemma}
\begin{proof}
First, uniform position principle \cite{ACGH} states that if $C_1\subseteq\P^a$ is an irreducible, nondegenerate curve then for a general hyperplane $H_a\subseteq\P^a$, the set $H_a\cap C_1$ spans $H_a$. Secondly, we want to show that $(*)$ forces these hyperplanes to satisfy $(ii)$ in the lemma. Choose a hyperplane $H_a$ which satisfies $(i)$ and suppose that the set $(H_a.C)$ generates a plane $\Pi_b\subseteq\P^b$ of codimension at least two. Hence for all hyperplanes $H_b$ passing through $\Pi_b$ we have $(H_a.C)\subseteq H_b$ and therefore $H_a\cap C_1\subseteq (C.H_b)$. Now $(*)$ says that the set $(C.H_b)$ lie in a hyperplane. Together with $(i)$ we get that $(C.H_b)\subseteq H_a$, for all hyperplanes $H_b$ containing $\Pi_b$. As the reunion of all hyperplanes containing $\Pi_b$ covers $\P^b$ we have that $C\subseteq H_a\times\P^b$, which is a contradiction. Lastly, it is easy to see that condition $(*)$ together with $(i)$ and $(ii)$ implies that condition $(iii)$ is also satisfied and we finish the proof.
\end{proof}
\begin{proof}[Proof of Theorem \ref{Theorem 1}]
The paragraph following Lemma \ref{Lemma2} says that we might not be able to obtain good projections only if the curve $C\subseteq\P^a\times\P^b$ satisfies property $(*)$. Now Lemma \ref{Lemma3} and $(*)$ implies that there exists an open set of hyperplanes $H_a\subseteq\P^a$ and another corresponding one of $H_b\subseteq\P^b$ such that condition $(iii)$ in Lemma \ref{Lemma3} is satisfied. Now choose $H_a$ and $H_b$ where $H_a\times\P^b$ and $\P^a\times H_b$ intersects transversally the curve $C$ at each point. We deduce, by condition $(iii)$ in Lemma \ref{Lemma3}, that $L_1\simeq L_2$. Denote this vector bundle by $L$. The openess of this condition says that there exists an open set of sections in $W_1$ and an open set of sections in $W_2$ which corresponds to each other. This forces $W_1=W_2$ inside $H^0(C,L)$, therefore $a=b$ and the curve $C$ is included in the graph of an automorphism and Theorem \ref{Theorem 1} is proved.
\end{proof}
\begin{remark} \label{Remark 2}
Choose a general codimension two plane $\Lambda_a$  and suppose it satisfies condition $(2_{a})$ from Lemma \ref{Lemma1}. As this condition is open then for a general codimension two plane in the first factor the same property is satisfied. Therefore for a curve $C\subseteq\P^a\times\P^b$ as in Theorem \ref{Theorem 1} we have plenty of pairs $(\Lambda_a ,\Lambda_b )$ which define a good projection for $C$.
\end{remark}

\section{Regularity bounds in general case}
In this section our goal is to prove that the bound on the multigraded regularity given in Theorem A holds for all curves $C\subseteq\P^a\times\P^b$, which satisfy the conditions in Theorem~\ref{Theorem 1}.
\begin{theorem}\label{Theorem 2}
Let $C\subseteq\P^a\times\P^b$ be a smooth curve of bidegree $(d_1,d_2)$ as in Theorem \ref{Theorem 1}, then the ideal sheaf $\mathcal{I}_{C|\P^a\times\P^b}$ is $(d_2-b+1,d_1-a+1)$-regular.
\end{theorem}
\noi The key to Theorem \ref{Theorem 2} is the following result, which will later allow us to connect the regularity of the ideal sheaf of $C$ with the regularity of a certain vector bundle on $\P^1\times\P^1$.
\begin{Proposition} \label{Proposition 1}
If the pair $(\Lambda_a,\Lambda_b)$ defines a good projection for $C$ then on $\P^1\times\P^1$ we have the following short exact sequence:
\begin{equation}
 0\rightarrow \mathcal{E} \rightarrow V_a\otimes \mathcal{O}_{\P^1\times\P^1}(-1,0)\oplus V_b\otimes \mathcal{O}_{\P^1\times\P^1}(0,-1)\oplus \mathcal{O}_{\P^1\times\P^1} \overset{\epsilon}\rightarrow (f_{a,b})_*(\mathcal{O}_C) \rightarrow 0
\end{equation}
where $\mathcal{E}$ is a vector bundle of rank $a+b-1$ and $V_a$, $V_b$ are vector spaces of dimension $a-1$ and $b-1$ respectively.
\end{Proposition}
\begin{proof} Blow-up $\P^a$ along $\Lambda_a$ and $\P^b$ along $\Lambda_b$ to get the diagram:
\[
\xymatrix{ & C\subseteq Y:=Bl_{\Lambda_a}(\P^a)\times Bl_{\Lambda_b}(\P^b) \ar[dl]_{\mu_a\times \mu_b} \ar[dr]^{p_a\times p_b} & \\
C\subseteq \P^a\times\P^b & & \bar{C}\subseteq \P^1\times\P^1\\}
\]
The morphism $p_a\times p_b$ will resolve the projection map $\pi_a\times\pi_b$, whose restriction to $C$ is a good projection for the curve. Now set:
\[
A_1:=(\mu_a\times\mu_b)^*(\mathcal{O}_{\P^a\times\P^b}(1,0))\text{ and }A_2:=(\mu_a\times\mu_b)^*(\mathcal{O}_{\P^a\times\P^b}(0,1))
\]
As $\Lambda_a\cap C_1=\{\varnothing\}$ and $\Lambda_b\cap C_2=\{\varnothing\}$, we can consider that $\mathcal{I}_{C|Y}=(\mu_a\times\mu_b)^*(\mathcal{I}_{C|\P^a\times\P^b})$. Using the diagram and the notations we made, we have two exact sequences on $\P^1\times\P^1$:
\[
0 \rightarrow  (p_a\times p_b)_*(\mathcal{I}_{C|Y}(A_i)) \rightarrow (p_a\times p_b)_*(\mathcal{O}_Y(A_i))\overset{\epsilon_i}\rightarrow (f_{a,b})_*(L_i)\text{ where } i=1,2
\]
We will be interested to find points $(x,y)\in\P^1\times\P^1$, where the stalk of either $\epsilon_1$ or $ \epsilon_2$ is surjective. In the case of $\epsilon_1$, by Nakayama's lemma, it suffices to show that the map: 
\[
\xymatrix{
(p_a\times p_b)_*(A_1)\otimes \mathbb{C}(x,y)  \ar[d]_{\simeq}  \ar[r]^{ \ \ \epsilon_1\otimes \mathbb{C}(x,y)} & (f_{a,b})_*(L_1)\otimes \mathbb{C}(x,y) \ar[d]^{\simeq} \\
H^0(\overline{\Lambda_a x}\times\overline{\Lambda_b y},\mathcal{O}_{\overline{\Lambda_a x}\times\overline{\Lambda_b y}}(1,0))  \ar[r] & 
H^0(\mathcal{O}_{C\cap \overline{\Lambda_a x}\times\overline{\Lambda_b y}}(1,0))
  }
\]
is surjective. Equivalently, we need to study the surjectiveness of the bottom horizontal map. We know that the pair $(\Lambda_a,\Lambda_b)$ defines a good projection, hence the intersection $C\cap \overline{\Lambda_ax}\times\overline{\Lambda_by}$ consists of at most two points. If it is a point then clearly the bottom horizontal map is surjective. If it consists of $(x^1,y^1)$ and $(x^2,y^2)$ then the bottom horizontal map is surjective whenever $x^1\neq x^2$. This implies that for all $(x,y)\in\P^1\times\P^1$ the stalk of at least one of the maps $\epsilon_1$ or $\epsilon_2$ is surjective. 

Note that $Bl_{\Lambda_a}(\P^a)=\P(V_a\otimes \mathcal{O}_{\P^1}\oplus \mathcal{O}_{\P^1}(1))$ and $Bl_{\Lambda_b}(\P^b)=\P(V_b\otimes \mathcal{O}_{\P^1}\oplus \mathcal{O}_{\P^1}(1))$, where $V_a$ and $V_b$ are vector spaces of dimension $a-1$ and $b-1$ respectively. This allows one to have the isomorphisms:
$$(p_a\times p_b)_*(A_1)\simeq V_a\otimes \mathcal{O}_{\P^1\times\P^1}\oplus \mathcal{O}_{\P^1\times\P^1}(1,0)\text{ and }(p_a\times p_b)_*(A_2)\simeq V_b\otimes \mathcal{O}_{\P^1\times\P^1}\oplus \mathcal{O}_{\P^1\times\P^1}(0,1)$$
Bear in mind the ideas from before, tensor $\epsilon_1$ by $\mathcal{O}_{\P^1\times\P^1}(-1,0)$, $\epsilon_2$ by $\mathcal{O}_{\P^1\times\P^1}(0,-1)$, and sum them together to get the following surjective map:
\begin{align}\notag
V_a\otimes \mathcal{O}_{\P^1\times\P^1}(-1,0)\oplus \mathcal{O}_{\P^1\times\P^1}\oplus V_b\otimes \mathcal{O}_{\P^1\times\P^1}(0,-1)\oplus \mathcal{O}_{\P^1\times\P^1} 
\overset{\epsilon_0}\longrightarrow (f_{a,b})_*(\mathcal{O}_C)
\end{align}
Notice that the second and the fourth component of the domain of $\epsilon_0$ have the same image, so we actually get the following short exact sequence:
\[
 0\rightarrow \mathcal{E} \rightarrow V_a\otimes \mathcal{O}_{\P^1\times\P^1}(-1,0)\oplus V_b\otimes \mathcal{O}_{\P^1\times\P^1}(0,-1)\oplus \mathcal{O}_{\P^1\times\P^1} \overset{\epsilon}\rightarrow (f_{a,b})_*(\mathcal{O}_C) \rightarrow 0
\]
with $\mathcal{E}=\text{Ker}(\epsilon)$. Now provided that $(f_{a,b})_*(\mathcal{O}_C)$ is Cohen-Macaulay sheaf with support of codimension 1 we have that $\mathcal{E}$ is a vector bundle of rank $a+b-1$ and this ends the proof.
\end{proof}
Now the idea is to find bounds on the multigraded regularity of the vector bundle $\mathcal{E}$, but before that we have the following proposition:
\begin{Proposition} \label{Proposition 2}
In the assumptions of Theorem \ref{Theorem 1}, we can make a choice of a good projection so that the dual vector bundle $\mathcal{E}^*$ is $(-1,-1)$-regular.
\end{Proposition}
\begin{proof}
Serre duality and the short exact sequence $(1)$ imply the vanishings: 
\[
H^2(\P^1\times\P^1,\mathcal{E}^*(-3,-1))=H^2(\P^1\times\P^1,\mathcal{E}^*(-1,-3))=0
\]
Again by Serre duality $H^2(\P^1\times\P^1,\mathcal{E}^*(-2,-2))=H^0(\P^1\times\P^1,\mathcal{E})$. Now use the sequence:
\[
 0\rightarrow H^0(\P^1\times\P^1,\mathcal{E}) \rightarrow H^0(V_a\otimes \mathcal{O}_{\P^1\times\P^1}(-1,0)\oplus V_b\otimes \mathcal{O}_{\P^1\times\P^1}(0,-1)\oplus \mathcal{O}_{\P^1\times\P^1}) \rightarrow H^0(C,\mathcal{O}_C) 
\]
to obtain that $H^0(\P^1\times\P^1,\mathcal{E})=H^1(\P^1\times\P^1,\mathcal{E})=0$.

It remains to show the vanishings of $H^1(\P^1\times\P^1,\mathcal{E}^*(-1,-2))$ and $H^1(\P^1\times\P^1,\mathcal{E}^*(-2,-1))$ respectively. We will prove that the first group is zero, as the second vanishing follows from the same ideas. First, by Serre duality it is isomorphic to $H^1(\P^1\times\P^1,\mathcal{E}(-1,0))$. Now use the vanishings above and the exact sequence:
\[
 0 \longrightarrow \mathcal{E}(-1,0) \longrightarrow \mathcal{E} \longrightarrow  \mathcal{E}|_{\{ x\}\times\P^1} \longrightarrow 0
\]
to get that $H^0(\{ x\}\times\P^1,\mathcal{E}|_{\{ x\}\times\P^1})=H^1(\P^1\times\P^1,\mathcal{E}(-1,0))$, $\forall x\in\P^1$. The morphism $p_1|_C$ has a nondegenerate image, so the multiplication map by $x$: $(f_{a,b})_*(\mathcal{O}_C)(-1,0) \longrightarrow (f_{a,b})_*(\mathcal{O}_C)$ is injective. Thus Snake lemma (Proposition 1.2, \cite{Bou}) implies that $\forall x\in\P^1$ we have the exact sequence:
\[
0\rightarrow \mathcal{E}|_{\{ x\}\times\P^1}\rightarrow V_a\otimes \mathcal{O}_{\{ x\}\times\P^1}\oplus V_b\otimes \mathcal{O}_{\{ x\}\times\P^1}(-1)\oplus \mathcal{O}_{\{ x\}\times\P^1}\rightarrow (f_{a,b})_*(\mathcal{O}_C)|_{\{ x\}\times\P^1}\rightarrow 0
\]
To end the proof is enough to show that for some $x\in\P^1$ the map:
\[
V_a\otimes H^0(\mathcal{O}_{\{ x\}\times\P^1})\oplus H^0(\mathcal{O}_{\{ x\}\times\P^1}) \overset{l_1\oplus l_2}\longrightarrow 
H^0((f_{a,b})_*(\mathcal{O}_C)|_{\{ x\}\times\P^1})
\]
is injective. For this purpose, suppose the projection map is given by the formula:
\[
 \pi_{a}\times\pi_{b}([x_0:...:x_a]\times [y_0:...:y_b])=[x_0:x_1]\times [y_0:y_1]
\] 
Remark \ref{Remark 2} implies that we can choose $x\in\P^1$ which satisfies:
\[
\overline{\Lambda_ax}\times\P^b\cap C=\{P_1,...,P_{d_1}\}\text{ with }P_i=(x^i,y^i)\text{ and }x^i\neq x^j\text{ for }i\neq j
\]
At last assume that $P_1\in\{x_2=...=x_a=0\}\times\P^b$. In this case we have:
\[
(f_{a,b})_*(\mathcal{O}_C)|_{\{ x\}\times\P^1}=\oplus_{i=1}^{d_1} \mathbb{C}_{P_i}
\]
and we can write $l_2(1)=(1,...,1) \text{ and } l_1(e_i)=(x_i|_{P_1},...,x_i|_{P_{d_1}})$ for a basis $\{ e_2,...,e_a\}$ of $V_a$. As $x_i|_{P_1}=0, \forall i=2...a$, it is enough to prove that $l_1$ is injective. Suppose the opposite, then there exists $(u_2,...,u_a)\in\mathbb{C}^{a-1}$ such that: 
\[
(\sum_{i=2}^au_ix_i|_{P_1},...,\sum_{i=2}^au_ix_i|_{P_{d_1}})=(0,...,0)
\]
This means that the set $\{x^1,...,x^{d_1}\}\subseteq\{ \sum_{i=2}^au_ix_i=0\}\cap \overline{x\Lambda_a}$, therefore the points $x^1,...x^{d_1}$ span a plane $\Pi_a$ of codimension at least two. Choose $(x,y)\in C$ with $x\notin \Pi_a$ and a hyperplane $H_a$ containing $\overline{\Lambda_ax}$. Now the intersection $H_a\times\P^b\cap C$ consists of at least $d_1+1$ points. This is a contradiction as we assumed that $p_1|_C$ has a nondegenerate image. Therefore $l_1$ is injective and the vector bundle $\mathcal{E}^*$ is $(-1,-1)$-regular.
\end{proof}
The last ingredient necessary for Theorem \ref{Theorem 2} is the following lemma, which will later allow us to connect the multigraded regularity of $\mathcal{E}$ to the one of $\mathcal{E}^*$:
\begin{lemma} \label{Lemma B}
Let $\mathcal{G}_1$ and $\mathcal{G}_2$ be two locally free sheaves on $\P^1\times\P^1$. If $\mathcal{G}_1$ is $(p,q)$-regular and $\mathcal{G}_2$ is $(m,n)$-regular then $\mathcal{G}_1\otimes \mathcal{G}_2$ is $(p+m,q+n)$-regular. In particular as we work over complex numbers then for all $k\in\mathbb{N}$ we have $\bigwedge^k (\mathcal{G}_1)$ is $(kp,kq)$-regular.
\end{lemma}
\begin{proof}
It is enough to consider the case when $\mathcal{G}_1$ and $\mathcal{G}_2$ are $(0,0)$-regular and we will only prove the vanishing of  $H^1(\P^1\times\P^1,\mathcal{G}_1\otimes \mathcal{G}_2(-1,0))$, as the other ones follow from the same ideas. Recall (Theorem 1.4 \cite{MaclaganSmith}) which implies that both $\mathcal{G}_1$ and $\mathcal{G}_2$ are globally generated, one therefore has the two short exact sequences:
\[
 0 \longrightarrow \mathcal{M}_i \longrightarrow \oplus \mathcal{O}_{\P^1\times\P^1} \longrightarrow \mathcal{G}_i \longrightarrow 0
\]
where $\mathcal{M}_i$ are locally free sheaves for $i=1,2$. Now tensor the second sequence with $\mathcal{G}_1(-1,0)$ and get the exact sequence in cohomology:
\[
 \oplus H^1(\P^1\times\P^1,\mathcal{G}_1(-1,0)) \rightarrow H^1(\P^1\times\P^1,\mathcal{G}_1\otimes \mathcal{G}_2(-1,0)) \rightarrow H^2(\P^1\times\P^1,\mathcal{G}_1\otimes \mathcal{M}_2(-1,0))
\]
As $\mathcal{G}_1$ is $(0,0)$-regular, the right group vanishes and is enough to prove that the left one also does. For this, tensor the first short exact sequence with $\mathcal{M}_2(-1,0)$ to get:
\[
 \oplus H^2(\P^1\times\P^1,\mathcal{M}_2(-1,0))\rightarrow H^2(\P^1\times\P^1,\mathcal{M}_2\otimes \mathcal{G}_1(-1,0)) \rightarrow 0
\]
and conclude that we only need $H^2(\P^1\times\P^1,\mathcal{M}_2(-1,0))=0$. Going back to the second short exact sequence, tensor it with $\mathcal{O}_{\P^1\times\P^1}(-1,0)$ and get the exact sequence:
\[
 H^1(\P^1\times\P^1,\mathcal{G}_2(-1,0))\rightarrow H^2(\P^1\times\P^1,\mathcal{M}_2(-1,0)) \rightarrow \oplus H^2(\mathcal{O}_{\P^1\times\P^1}(-1,0))
\]
As $\mathcal{G}_2$ is $(0,0)$-regular this implies our vanishing, and the proof is done.
\end{proof}

\begin{proof}[Proof of Theorem \ref{Theorem 2}]
Theorem \ref{Theorem 1} and Remark \ref{Remark 2} implies the existence of "plenty" of good projections. This instead by Proposition \ref{Proposition 1} helps us to construct a vector bundle $\mathcal{E}$, which dual $\mathcal{E}^*$ is $(-1,-1)$-regular by Proposition \ref{Proposition 2}. Using the isomorphism:
\[
\mathcal{E}\simeq \bigwedge^{a+b-2}(\mathcal{E}^*)\otimes \text{det}(\mathcal{E})
\]
with Lemma \ref{Lemma B} and sequence $(1)$ in Proposition \ref{Proposition 1}  we obtain that $\mathcal{E}$ is $(d_2-b+1,d_1-a+1)$-regular. If we write $l:=d_2-b$ and $n:=d_1-b+1$ this says that the map:
\begin{align}\notag
H^0(V_1\otimes \mathcal{O}_{\P^1\times\P^1}(l-1,k)\oplus V_2\otimes \mathcal{O}_{\P^1\times\P^1}(l,k-1)\oplus \mathcal{O}_{\P^1\times\P^1}(l,k)) \overset{g}\longrightarrow H^0(\mathcal{O}_C(l,k)) 
\end{align}
is surjective. Assume again that the projection map is given by the formula:
\[
(\pi_{\Lambda}\times \pi_{\Gamma})([x_0:...:x_a]\times [y_0:...:y_b])=[x_0:x_1]\times [y_0:y_1]
\]
Thus the fact that $g$ is surjective, says that $H^0(C, O_C(l,k))$ is generated by the restriction to $C$ of polynomials of the type $x_iF$, $y_jG$ and $H$, $i=2...a$, $j=2...b$, where $F\in H^0(\mathcal{O}_{\P^1\times\P^1}(l-1,k))$, $G\in H^0(\mathcal{O}_{\P^1\times\P^1}(l,k-1))$ and $H\in H^0(\mathcal{O}_{\P^1\times\P^1}(l,k))$. Hence the map:
\[
H^0(\P^a\times\P^b, \mathcal{O}_{\P^a\times\P^b}(d_2-b,d_1-a+1))\rightarrow H^0(C,\mathcal{O}_C(d_2-b,d_1-a+1)) 
\]
is surjective and we have $H^1(\P^a\times\P^b,\mathcal{I}_{C|\P^a\times\P^b}(d_2-b,d_1-a+1))=0$. Symmetrically we have the vanishing of $H^1(\P^a\times\P^b,\mathcal{I}_{C|\P^a\times\P^b}(d_2-b+1,d_1-a))$. Now let's use the sequence:
\[0\longrightarrow \mathcal{I}_{C|\P^a\times\P^b}\longrightarrow \mathcal{O}_{\P^a\times\P^b}\longrightarrow \mathcal{O}_C\longrightarrow 0
\]
together with sequence $(1)$ from Proposition \ref{Proposition 1} to get:
\[
H^2(\P^a\times\P^b,\mathcal{I}_{C|\P^a\times\P^b}(m,n))=H^1(C,\mathcal{O}_C(m,n))=H^2(\P^1\times\P^1,\mathcal{E}(m,n))\text{   } \forall m,l\geq -1
\]
This implies the vanishings of $H^2$-s. To finish the proof notice that $d_2-b+1\geq 1$ and $d_1-a+1\geq 1$. Using vanishings of line bundles on $\P^a\times\P^b$ we conclude that $\mathcal{I}_{C|\P^a\times\P^b}$ is $(d_2-b+1,d_1-a+1)$-regular.
\end{proof}

\section{Regularity bounds in special case}
In order to finish the proof of Theorem A, it remains to study the case when $r:=a=b$ and $C$ is included in the graph of an automorphism of $\P^r$. This is discussed in the following proposition:
\begin{Proposition}\label{Proposition 3}
Let $C\subseteq\P^r\times\P^r$ be a smooth curve of bidegree $(d,d)$. If $C$ is contained in the diagonal $\Delta_{\P^r}$ and is nondegenerate, then $\mathcal{I}_{C|\P^r\times\P^r}$ is $(d-r+1,d-r+1)$-regular.
\end{Proposition}
\begin{proof}
Denote by $J:=\mathcal{I}_{C|\P^r\times\P^r}$ and $k:=d-r+1$. As $C$ is included in the diagonal we have the isomorphism $O_{\P^r\times\P^r}(1,0)|_C\simeq \mathcal{O}_{\P^r\times\P^r}(0,1)|_C$. Call $L$ this line bundle of degree $d$. To prove that $H^2(\P^r\times\P^r,J(k-2,k))=0$ use the following short exact sequence:
$$0\longrightarrow J \longrightarrow \mathcal{O}_{\P^r\times\P^r} \longrightarrow \mathcal{O}_C \longrightarrow 0$$
If we twist the sequence by $\mathcal{O}_{\P^r\times\P^r}(k-2,k)$, then we have in cohomology the sequence:
$$H^1(\P^r\times\P^r,\mathcal{O}(k-2,k))\rightarrow H^1(C,L^{\otimes 2k-2})\rightarrow H^2(\P^r\times\P^r,J(k-2,k))\rightarrow H^2(\P^r\times\P^r,\mathcal{O}(k-2,k))$$
As the curve $C\subseteq\P^r$ is nondegenerate we have $k\geq 1$. Therefore the first and the last group vanish and the middle ones are isomorphic. Recall (\cite{GLP}) that for a nondegenerate curve $C\subseteq\P^r$ of degree $d$ we have $L^{\otimes n}$ is non-special for all $n\geq d-r$, where $L:=\mathcal{O}_{\P^r}(1)|_C$. In our case $2k-2=2d-2r\geq d-r$ and the vanishing follows. The same ideas are used to prove: 
\[
H^2(\P^r\times\P^r, J(k-1,k-1))=H^2(\P^r\times\P^r,J(k,k-2))=0
\]
It remains to show $H^1(\P^r\times\P^r,J(k-1,k))=0$. This vanishing is equivalent to the surjectiveness of the following map:
\[
H^0(\mathcal{O}_{\P^r}(k-1))\otimes H^0(\mathcal{O}_{\P^r}(k)) \rightarrow H^0(C,L^{\otimes 2k-1})
\]
Now this map can be factored as follows:
\[
\xymatrix{  H^0(\mathcal{O}_{\P^r}(k-1))\otimes H^0(\mathcal{O}_{\P^r}(k)) \ar[r] \ar[d]^{l} &  H^0(C,L^{\otimes 2k-1}) & \\
  H^0(\P^r,\mathcal{O}_{\P^r}(2k-1))\ar[ur]^{u}  & \\}
\]
As $l$ is surjective we only need to prove that $u$ is surjective. For this we use a result from \cite{GLP}, which states that for a nondegenerate curve $C\subseteq\P^r$ of degree $d$ we have: 
\[
H^0(\P^r,\mathcal{O}_{\P^r}(n))\rightarrow H^0(C,L^{\otimes n}) \text{ is surjective } \forall n\geq d-r+1
\]
In our case $n=2k-1=2d-2r+1\geq d-r+1$ and this finishes the proof.
\end{proof}
\section{Examples}
The paper ends with an example which shows that the bound we have in Theorem A is the best possible. The curve given in this example has the property that the sets reg$(C)$ and $(d_2-b+1,d_1-a+1)+\mathbb{N}^2$ coincide. Thus by Theorem A, $(d_2-b+1,d_1-a+1)+\mathbb{N}^2$ is the maximal set contained in $\text{reg}(C)$ for all curves $C\subseteq \P^a\times\P^b$ of bidegree$(d_1,d_2)$ with nondegenerate birational projections.

The idea is to find examples of curves with high order "secant lines". In our case we will consider the "secant lines" of the following type:
\[
l\times [y_0:...:y_b]\subseteq\P^a\times\P^b, \text{ where }l\subseteq\P^a \text{ is a line }
\]
Suppose that $l\times[y_0:...:y_b]\cap C$ consists of $k$ points. If $s\in H^0(\P^a\times\P^b,O_{\P^a\times\P^b}(k-1,n))$ is a hypersurface, vanishing along the curve $C$, then the index of intersection with the "secant line" $l\times[y_0:...:y_b]$ is $k-1$. Thus $s$ also vanishes along $l\times[y_0:...:y_b]$. It follows that the sheaf $\mathcal{I}_{C|\P^a\times\P^b}(k-1,n)$ is not globally generated $\forall n\in\mathbb{N}$, and therefore (Theorem 1.4 \cite{MaclaganSmith}) it is not $(k-1,n)$-regular. With this in hand we have the following example:
\begin{example} Let the morphism $\psi :C=\P^1\longrightarrow \P^2\times\P^2$ be given by the formula:
\[
\psi ([s:t])=[t^2s-4s^3:t^3-4s^2t:t^2s-3s^3]\times [s^2t-t^3:s^3-st^2:t^3]
\]
First notice that $\psi$ defines an embedding, such that the curve $C$ is of bidegree $(3,3)$ with nondegenerate birational projections. At the same time we have the following:
\[
C\cap \{x_2=0\}\times [0:0:1]=\{[1:1],[1:-1]\}
\]
\[
C\cap [0:0:1]\times\{4y_0+3y_2=0\}=\{[1:2],[1:-2]\}
\]
Bearing in mind the ideas above we obtain that the ideal sheaf $\mathcal{I}_{C|\P^2\times\P^2}$ is not $(1,s)$ and $(t,1)$-regular for all $s,t\in \mathbb{N}$. As Theorem A states that this ideal is $(2,2)$-regular, we conclude that $\text{reg}(C)=(2,2)+\mathbb{N}^2$. It is easy to notice that we can generalize this example i.e. find a rational curve $C\subseteq\P^r\times\P^r$ of bidegree $(r+1,r+1)$ and nondegenerate birational projections such that $\text{reg}(C)=(2,2)+\mathbb{N}^2$. 
\end{example}

\end{document}